\documentclass{IEEEtran}
        
\usepackage{hyperref}
\usepackage{graphicx}
\usepackage{amssymb}
\usepackage{epstopdf}
\usepackage{amsthm}
\usepackage{amsmath}
\usepackage[dvipsnames]{xcolor}
\usepackage{algorithm}
\usepackage[utf8]{inputenc}
\usepackage[english]{babel}
\usepackage{algpseudocode}
\newtheorem{theorem}{Theorem}
\newtheorem{lemma}{Lemma}

\graphicspath{{Figures/}}

\title{On the Capacity of the Quantum Switch \\ with and without Entanglement Decoherence}
\author{
    \IEEEauthorblockN{V\'ictor Valls\IEEEauthorrefmark{1}, Panagiotis Promponas\IEEEauthorrefmark{2}, Leandros Tassiulas\IEEEauthorrefmark{2}} \\
    \IEEEauthorblockA{\IEEEauthorrefmark{1}IBM Research Europe -- Dublin
    }
    \IEEEauthorblockA{\IEEEauthorrefmark{2}Department of Electrical Engineering, Yale University
    \\}
}

\begin{document}

\maketitle

\begin{abstract}
This paper studies the capacity of the quantum switch for two decoherence models: when link-level entanglements last  (i) for a time slot, or (ii) until they are used to serve a request (i.e., there is no decoherence). The two models are important as they set lower and upper bounds on the capacity region for any other decoherence model. The paper's contributions are to characterize the switch capacity region for both decoherence models and to propose throughput-optimal policies based on gradient descent. 
\end{abstract}


\section{Introduction}

\begin{figure}
\centering
\includegraphics[width=0.95\columnwidth]{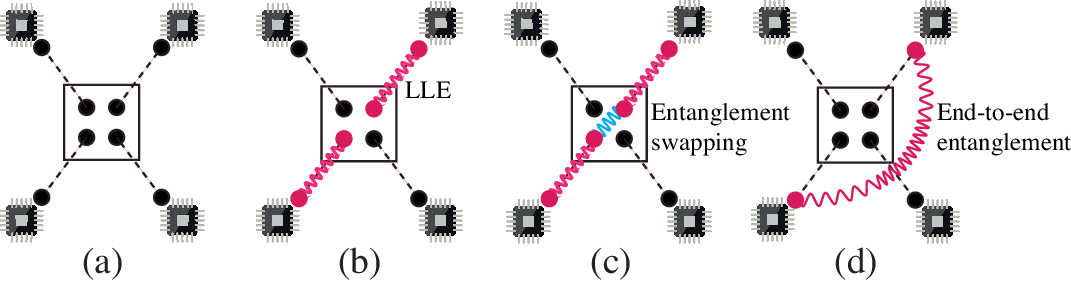}
\caption{Illustrating the operation of a quantum switch. LLEs are created between the clients (e.g., quantum processors) and the switch (b), and then the switch performs an entanglement swapping (c) to create an end-to-end entanglement (d).}
\label{fig:quantum_switch_big_picture}
\end{figure}

We study the capacity of the quantum switch and the design of throughput-optimal policies. In brief, a quantum switch is a star graph connected to clients with dedicated links (see Figure \ref{fig:quantum_switch_big_picture}) \cite{VGNT21}. The links are used for creating link-level entanglements (LLEs) between the switch and the clients (Figure \ref{fig:quantum_switch_big_picture}b), and the task of the switch is to perform entanglement swappings to create end-to-end entanglements (Figures \ref{fig:quantum_switch_big_picture}c \& \ref{fig:quantum_switch_big_picture}d). End-to-end entanglements are used by quantum applications to, for example, teleport information qubits.\footnote{The quantum switch in this paper is different from another device also referred as quantum switch, which controls the order in which quantum messages traverse communication channels. See, for example,  \cite{CC20}.} 

The problem we want to solve is the following. For a given arrival process of end-to-end entanglement requests, we want to find an entanglement swapping policy that maximizes the switch's utilization. This problem is different from scheduling problems in wired and wireless networks since LLEs are generated randomly over time and consumed to serve requests. 
Moreover, LLEs can have different durations depending on the underlying technology, which complicates the characterization of the switch capacity region (i.e., the set of loads that can be supported) and the design of throughput-optimal policies.

Previous work has studied the operation of a quantum switch under various settings. The papers in \cite{VGNT21} and \cite{NVGT22} investigate an idealized switch with bipartite and tripartite end-to-end entanglement requests with symmetric loads and negligible decoherence \cite{NC10}. The work in \cite{DRT21} studies a quantum switch with bipartite requests when there is no memory decoherence and LLE attempts succeed probabilistically. In particular,  \cite{DRT21} characterizes the switch capacity region when the switch has unlimited memory and studies the impact of finite memory numerically. Finally,  \cite{VT22} analyzes the capacity of the quantum switch when LLEs decohere after a time slot, which is similar to wireless network models with time-varying connectivity \cite{GNT06}.

In this paper, we study the capacity of the quantum switch for two decoherence models: 
when LLEs last (i) for a time slot (Sec.\ \ref{sec:maxdecoherence}) or (ii) until they are used to serve a request, i.e., there is no decoherence (Sec.\ \ref{sec:model_without_decoherence}).
Studying these models is important because they correspond, respectively, to the smallest and the biggest set of arrival loads that a quantum switch can support. The paper's contributions are to characterize the switch capacity region for both decoherence models and to propose throughput-optimal policies based on gradient descent (Algorithms \ref{al:spdga} and \ref{al:cc_spdga}). Our results differ from previous work because  our model and algorithms encode the memory constraints of practical systems. In particular, when the LLEs last for a time slot, the switch memory has to be equal to the number of LLEs generated during that time. And when there is no decoherence, the memory required is a system parameter that affects the algorithm's optimality, i.e., the maximum load that the switch can support.  Finally,  our technical approach is based on convex optimization and gradient descent for \emph{both} models, which underlines the fundamental differences between the models and the technical challenges. 

\section{Decoherence Model: \\LLEs Last For One Time Slot}
\label{sec:maxdecoherence}

\subsection{Quantum switch model and operation}
\label{sec:model} 

\subsubsection{Model}
We present a system model in line with the models in \cite{DRT21, VT22, PVTT22}. In brief, consider a quantum switch connected to $N$ clients with dedicated links. The switch operates in slotted time, where, in each time slot $k=1,2,3,\dots$,   
\begin{itemize}
\item [(i)] clients attempt to create LLEs with the switch; and
\item [(ii)] the switch performs an entanglement swapping with the clients' LLEs. 
\end{itemize}
Without loss of generality, we assume that clients only attempt one LLE per time slot, and that these succeed with probability $\tau_j$, $j \in \{1,\dots,N\}$. The successful LLEs are stored in separate quantum memories; one for each client. In this section, we assume that LLEs last for \emph{one time slot} due to decoherence, where the duration of a time slot is equal to the time required to carry out (i) and (ii).  For example, if it takes 1 unit of time to generate an LLE, and 2 units to perform an entanglement swapping, the duration of a time slot is 3. Also, we assume that LLEs are generated \emph{only} at the beginning of a time slot, and that if an LLE has not been used when the time slot ends, this is discarded. Hence, the switch only needs to store one LLE per client at a time.

Requests for end-to-end entanglements connect pairs of clients, and they arrive at the switch following a stochastic process $\{b_k\}_{k=1}^\infty$, where $b_k \in \mathbf N^d$ is bounded with 
\begin{align}
 b:= \lim_{k \to \infty} \frac{1}{k} \sum_{i=1}^k b_i, \label{eq:average_b}
\end{align}
and $d := N (N-1)/2$. Also, we assume $\{b_k\}_{k=1}^\infty$ is i.i.d.\ and independent of the clients' successful LLE attempt probabilities.
The requests are stored in queues $Q_k$, which evolve as follows 
\begin{align}
Q_{k+1} = [Q_k + b_k - y_k]^+.
\label{eq:queue_evolution}
\end{align}
We use vector $y_k \in Y_k \subseteq \{0,1\}^d$ to denote the service of requests at time $k$ obtained as a result of performing an entanglement swapping (see Figure \ref{fig:quantum_switch_big_picture}). 
The set 
\begin{align}
Y_k \subseteq \{0,1\}^d \label{eq:action_set}
\end{align}
contains the admissible service vectors (or, actions) in a time slot since  it is not possible to serve a request unless all the clients that participate in it have an active LLE (see Figure \ref{fig:quantum_switch_big_picture}). Also, note that since there is only one LLE per client in a time slot, the switch can serve multiple requests in a time slot but only one request per client. Finally, an entanglement swapping consumes the clients' LLEs, and so an LLE can only be used to serve one request. 

\subsubsection{Operation}
The switch's task is to carry out entanglements swappings (i.e., select $y_k \in Y_k$) to serve as many end-to-end entanglement requests as possible.
This type of control problem is often tackled with Lyapunov optimization techniques where the goal is to stabilize a system of queues (e.g., \cite{GNT06}), i.e., $\lim_{k \to \infty} \frac{1}{k} \sum_{i=1}^k \mathbf E[Q_i] \prec \infty$. In this paper, we use an equivalent approach based on convex optimization and  gradient descent  (e.g., \cite{CLC06, VL19}) where 
\begin{itemize}
\item[(i)] Lagrange multipliers play the role of queues; and 
\item [(ii)] Stability will be a consequence of satisfying the long-term or ``average'' constraint:
\end{itemize}
\begin{align}
b : = \lim_{k \to \infty} \frac{1}{k} \sum_{i=1}^k b_i \preceq  y:= \lim_{k \to \infty} \frac{1}{k} \sum_{i=1}^k y_i. \label{eq:long-term-constraint}
\end{align}

The approach consists of two steps. We first (i) formulate the fluid or static problem that is decoupled from the switch operation (Sec.\ \ref{sec:fluid_formulation_with_decoherence}), and then (ii) design a numerical method to solve the static problem that can be used to operate the quantum switch  (Sec.\ \ref{sec:alg_with_decoherence}). In particular, an iteration in the numerical method corresponds to a time slot in the quantum switch, and the capacity region of the switch is the set of feasible service vectors (i.e., solutions) to the fluid problem.

\subsection{Static formulation}
\label{sec:fluid_formulation_with_decoherence}

\subsubsection{Set of service vectors (capacity region)} We proceed to characterize the set of average service rates that the switch can generate. 
In short, let $\mathbb S := \{0,1\}^N$ denote the switch states or LLEs connectivity with the clients in a time slot. In particular, a $1$ in the $j$'th entry of a state vector $s \in \mathbb S$ represents that the $j$'th client has an active LLE with the switch, and 0 otherwise. Each state $s \in \mathbb S$ is associated with a probability $p(s)$ that depends on the client's successful LLE attempt probabilities $\tau_j$, $j\in \{1,\dots,N\}$. For example, the switch will be connected to all clients in a time slot with probability $p(s = \mathbf 1) = \prod_{j=1}^N \tau_j$.

Similarly, let $Y(s) \subseteq \{0,1\}^d$ denote the set of entanglement swapping operations available in a state $s$. 
The average service rate in a state $s$ is given by 
 \begin{align}
 x(s) := \sum_{y \in Y(s)} \theta(y, s) y, \label{eq:xs}
 \end{align}
where $\sum_{y \in Y(s)} \theta(y, s) = 1$, $\theta(y, s) \ge 0$ for all $y \in Y(s)$. That is, $x(s)$ is a weighted average of the feasible service vectors in state $s$. The \emph{total} service rate of the switch is the weighted average of the rates in each of the states. Let  
\[
X(s) := \mathrm{conv}(Y(s))
\]
be the convex hull of $Y(s)$ containing all possible vectors $x(s)$ as defined in \eqref{eq:xs}.
Then, the set of all feasible service rate vectors (i.e., the switch capacity region) is
\begin{align}
C:=  \left\{ u \in \mathbf R^d \mid u =   \sum_{s \in \mathbb S} p(s) x(s), x(s) \in X(s) \right\}.
\label{eq:capacity_region}
\end{align}

\subsubsection{Optimization problem} The static or fluid problem we want to solve is the following. For an average arrival vector $b$ and some fixed system states probabilities $p(s)$, the objective is to find the service rates vectors in each state (i.e., $x(s)$)  to serve all the requests. That is, solve:
\begin{align}
\begin{tabular}{ll}
$\underset{x(s) \in X(s), \! \ s \in \mathbb S}{\text{minimize}}$ & $1$\\
subject to & $\displaystyle b \preceq x:= \sum_{s \in \mathbb S} p(s) x(s)$
\end{tabular}
\label{eq:primal_problem}
\end{align}

In the next section, we present an algorithm that solves \eqref{eq:primal_problem} and can also operate the quantum switch with the queue dynamics in \eqref{eq:queue_evolution}. In particular, an iteration in the numerical method corresponds to a switch time slot $k$, and each time slot is associated with a system state $s_k \in \mathbb S$ and a request arrival vector $b_k$. 

\subsection{Algorithm}
\label{sec:alg_with_decoherence}

In short, we tackle \eqref{eq:primal_problem} by formulating its Lagrange dual problem \cite[Ch.\ 5]{BV04}:
\begin{align}
\begin{tabular}{ll}
$\underset{\lambda \succeq 0}{\text{maximize}}$ & $ h(\lambda) $ \\
\end{tabular}
\label{eq:lagrange_dual_problem}
\end{align}
where $
h(\lambda) : = 1+ 
 \left \langle \lambda, b  \right\rangle -  \sum_{s \in \mathbb S}  \min_{u \in X(s)} \left\langle  \lambda, p(s) u\right \rangle$ and $\lambda \in \mathbf R^d_+$ is a Lagrange multiplier. 
It is well-known that solving \eqref{eq:lagrange_dual_problem} is equivalent to solving \eqref{eq:primal_problem} when strong duality holds \cite[Ch.\ 5]{BV04}. And since the dual problem is a concave maximization problem, we can solve \eqref{eq:lagrange_dual_problem}  numerically with, for example, the projected (sub)gradient method (Algorithm \ref{al:pdga}).

We can use Algorithm \ref{al:pdga} to operate a quantum switch by making the (sub)gradient updates stochastic (Algorithm \ref{al:spdga}). In brief, a (sub)gradient of the dual function $h$ w.r.t.\ $\lambda_k$ is
\begin{align}
\nabla h(\lambda_k) = b  - \sum_{s \in \mathbb S} p(s) u(s,\lambda_k)
\end{align}
where $u(s, \lambda_k) \in \arg \min_{u \in X(s)} \langle \lambda_k, u \rangle$ for all $s \in \mathbb S$.
Next, observe that if $\{b_k\}_{k=1}^\infty$ and $\{ s_k \in \mathbb S \}_{k=1}^\infty$ are i.i.d.\ processes, then $b_k - y_k $ with $ y_k \in \arg \min_{v \in Y_k} -\langle \lambda_k, v \rangle$ is a stochastic (sub)gradient of $h$ at $\lambda_k$,\footnote{An extreme point is always a solution to a linear program on a convex polytope. Formally,  $\arg \min_{u \in X(s)} \langle -\lambda_k , u \rangle \supseteq \arg \min_{u \in Y(s)} \langle -\lambda_k , u \rangle \ne \emptyset$ for any $s \in \mathbb S$. See, for example, \cite{Dan63}.}
i.e., 
\begin{align}
\nabla h(\lambda_k) = \mathbf E_b \mathbf E_s [b_k - y_k].
\end{align}
Thus, the dual variables update in Algorithm \ref{al:spdga} with $\alpha  =1$ (c.f.\ step 4 in Algorithm \ref{al:pdga}) is
\begin{align}
\lambda_{k+1} =  [ \lambda_k + b_k - y_k ]^+, \label{eq:lupdate7}
\end{align}
which is equivalent to the queues dynamics in Eq.\ \eqref{eq:queue_evolution} if we identify $\lambda_k$ with $Q_k$. 
That is, we are solving the static problem in \eqref{eq:primal_problem} while operating the quantum switch. Also, note that the update in \eqref{eq:lupdate7} only needs to know the arrivals  (i.e., $b_k$) and the set of admissible actions (i.e., $Y_k$) in a time slot. 

The convergence of  Algorithm \ref{al:spdga} is standard from (sub)gradient methods (see, for example, \cite{VL19, NO09}). In particular, when the expected value of the requests arrivals (i.e., $b$) is in the interior of the capacity region $C$, Slater's condition \cite{BV04} is satisfied and Algorithm \ref{al:spdga} ensures that
\[
 \mathbf E \left[ \left\| \left[ \frac{1}{k} \sum_{i=1}^k b_i - \frac{1}{k} \sum_{i=1}^k y_i \right]^+ \right\|^2 \right] \le O\left(\frac{1}{\sqrt k}\right),
\] which means that the constraint in \eqref{eq:primal_problem} is satisfied asymptotically (see, for example, \cite[Thm.\ 1; claim (iii)]{VL19}). Furthermore,  Algorithm \ref{al:spdga} ensures that $\lim_{k \to \infty} \frac{1}{k} \sum_{i=1}^k \mathbf E [ \lambda_i ] $  is bounded, which is equivalent to having stable queues. 

\begin{algorithm}[t]
\caption{Projected Dual Gradient Ascent (PDGA)}
\begin{algorithmic}[1]
\State \textbf{Set:} $\lambda_0 = 0$, $\alpha > 0$
\While{stopping criterion is not met}
\State $x_k(s) \in \arg \min_{u \in X(s)} \ - \langle \lambda_k , u \rangle$ \quad  $\forall s \in \mathbb S$
\State $\lambda_{k+1} = \left[ \lambda_k +  \alpha \nabla h(\lambda_k) \right]^+$
\State $k \leftarrow k + 1$
\EndWhile
\end{algorithmic}
\label{al:pdga}
\end{algorithm}

\begin{algorithm}[t]
\caption{Stochastic PDGA for \eqref{eq:lagrange_dual_problem}}
\begin{algorithmic}[1]
\State \textbf{Set:} $\lambda_0 = 0$, $\alpha = 1$
\While{stopping criterion is not met}
\State Observe state $s_k$ and learn action set $Y_k = Y(s_k)$
\State $y_k \in \arg \max_{v \in {Y_k}} \ \langle \lambda_k, v \rangle$
\State $\lambda_{k+1} = \left[  \lambda_k +  b_k - y_k \right]^+$
\State $k \leftarrow k + 1$
\EndWhile
\end{algorithmic}
\label{al:spdga}
\end{algorithm}



The assumption that the system states and request arrivals are i.i.d. processes is essential. If LLEs lasted for multiple time slots, the system states $s_k \in \mathbb S$ would become coupled with the control policy, and so we cannot ensure $\{ s_k \}_{k=1}^\infty$ is an i.i.d.\ process, and as a result, guarantee that Algorithm \ref{al:spdga} solves the problem in \eqref{eq:primal_problem}. The latter is the technical reason why we model the problem differently when there is no decoherence.

\section{Decoherence Model: \\LLEs Last Until Used to Serve a Request}
\label{sec:model_without_decoherence}

In this section, we consider the case where there is no decoherence: LLEs last until they are used to serve a request. We cannot use the same model as in the previous section because when the LLEs do not expire after a time slot, the actions made by the network controller affect the future system states (i.e., the LLEs that will be available). As a result, \emph{the states and action sets are not i.i.d.}, which is a crucial property to ensure that the stochastic gradient method in Sec.\ \ref{sec:alg_with_decoherence} converges.

To tackle this case, we extend the model in the previous section (Sec.\ \ref{sec:extended-model}) and introduce a congestion control mechanism (Sec.\ \ref{sec:congestion-control}). The congestion mechanism will allow us to obtain a gradient algorithm that waits to serve a request until there are enough LLEs available (Sec.\ \ref{sec:cc_algorithm}).

\subsection{Model, capacity region, and fluid problem}
\label{sec:extended-model}

We extend the model in Sec.\ \ref{sec:model} to treat LLEs as a type of commodity. In short, let $\{ \hat b_k \}_{k=1}^\infty$ be the LLEs arrival process with $\hat b_k \in \mathbb \{0,1\}^N$, which are stored in $N$ separate queues $\hat Q_k$ in addition to the $d$ queues for requests. That is, there is a total $d + N$ queues. The queues evolve as follows $\hat Q_{k+1} = [\hat Q_k + \hat b_k - A y_k]^+$ where $y_k \in Y \subset \{ 0,1 \}^d$ indicates which requests to serve and matrix $A \in \{0,1\}^{N\times d}$ maps the LLEs required to serve a request.\footnote{The $j$'th column of $A$ is a binary vector where the entries with a $1$ indicate the LLEs required to serve the request associated with that column.} Importantly, now the action set $Y$ is the same for all time slots, and we assume that it is possible to serve multiple requests but only one request per client---as in Sec.\ \ref{sec:model}.

The capacity region in this setting is 
\begin{align}
C:= \{ x \in X \mid  Ax \preceq \hat b\},
\end{align}
where $X:=\mathrm{conv}(Y)$ and $\hat b : = \lim_{k \to \infty} \frac{1}{k} \sum_{i=1}^k \hat b_i$. In words, set $C$ contains the service vectors such that there are enough LLEs on \emph{average} to serve the requests. Thus, the static or fluid optimization problem is:
\begin{align}
\begin{tabular}{lllll}
$\underset{x \in X }{\text{minimize}}$ & $1$ \\
subject to & $ b \preceq x $, \ $Ax \preceq \hat b$
\end{tabular}
\label{eq:primal_problem_no_decoherence}
\end{align}
Note that \eqref{eq:primal_problem_no_decoherence} is similar to \eqref{eq:primal_problem}, but now the action set $Y$ does not depend on the system states and there is the additional constraint $Ax \preceq \hat b$. 

Our goal is to stabilize the queues $Q_k$ associated with the end-to-end entanglement requests.
However, we cannot use the approach in the previous section directly because---unlike \eqref{eq:long-term-constraint}---it is not enough to satisfy the constraint $Ax \preceq \hat b$ on ``average'' or in the long-term. The constraint must be satisfied in all time slots since we cannot serve a request if there are not enough LLEs available in the queues $\hat Q_k$. Such issue did not arise in Sec.\ \ref{sec:maxdecoherence} because the action sets \eqref{eq:action_set} contained only service vectors for which there are enough LLEs. To get around that issue, we implement a congestion control mechanism to make the gradient algorithm wait to serve a request until there are enough LLEs available.

\subsection{Congestion control mechanism}
\label{sec:congestion-control}
\begin{figure}
\centering
\includegraphics[width=0.95\columnwidth]{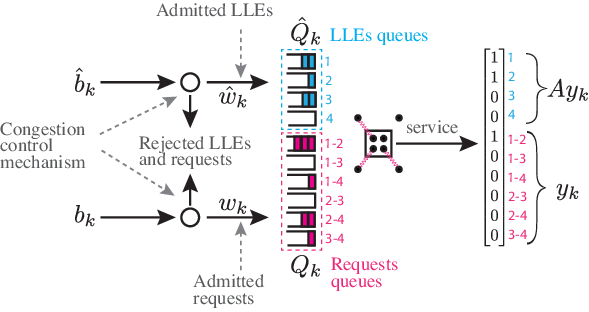}
\caption{Illustrating the congestion control model in Sec.\ \ref{sec:model_without_decoherence}. The congestion mechanism makes the gradient algorithm wait to serve a request until there are enough LLEs are available. For example, the algorithm will not serve a request from queue 1-4 (red) unless queues 1 and 4 (blue) have an LLE.}
\label{fig:congestion_control}
\end{figure}

The switch operation with the congestion control mechanism is shown in Figure \ref{fig:congestion_control}. Requests and LLEs arrive in each time slot, and the switch has to decide which to admit into the $d + N$ queues. The requests ($w_k \in \{0,1\}^d$) and LLEs ($\hat w_k \in \{0,1\}^N$) admitted must always be smaller than the arrivals (i.e., $w_k \preceq b_k$ and $\hat w_k \preceq \hat b_k$). The requests and LLEs  that are not admitted are discarded and cannot be readmitted in future time slots. 

As in Sec.\ \ref{sec:fluid_formulation_with_decoherence}, we use system states to model that the admission decisions depend on the requests and LLEs arrivals. 
Let $\mathbb S$ be a collection of system states and $b(s) \in \{0,1\}^d$ and $\hat b(s) \in \{0,1\}^N$ the requests and LLEs arrivals in a state $s \in \mathbb S$. The admission control variables in a state $s \in \mathbb S$ are
\begin{align}
B(s) & := \{ u \in \{0,1\}^d \mid u^{(j)} \in \{0, b^{(j)}(s)\} \} \\
\hat B(s) & := \{ \hat u \in \{0,1\}^N \mid \hat u^{(j)} \in \{0, b^{(j)}(s)\} \} 
\end{align}
where $u^{(j)}$ is the $j$'th component of vector $u$. Also, define $
Z:= \{z \in \mathbf R_+^d \mid z \preceq b \}$ and $\hat Z:= \{\hat z \in \mathbf R_+^d \mid \hat z \preceq \hat b \}$ be the sets with the ``average'' admission of requests and LLEs.\footnote{Set $Z$ can be equivalently defined as  $Z  := \textstyle \{  z \in \mathbf R^d \mid  z:= \sum_{s \in \mathbb S} p(s)  u(s),  u(s) \in \mathrm{conv}(B(s)) \}$. The same applies to set $\hat Z$. } The static or fluid congestion control problem is:
\begin{align}
\begin{tabular}{ll}
$\underset{x \in X,   z \in Z, \hat z \in \hat Z }{\text{minimize}}$ & $ \delta \underbrace{  \left\langle \mathbf 1, \begin{bmatrix} x \\ Ax \end{bmatrix} \right\rangle }_{\text{service}} -  \gamma \underbrace{ \left\langle  \mathbf 1, \begin{bmatrix}  z \\ \hat z \end{bmatrix} \right\rangle}_{\text{admission}}$ \\
subject to & $ \begin{bmatrix} z \\ \hat z \end{bmatrix}  \preceq \begin{bmatrix} x \\ Ax \end{bmatrix} \begin{matrix}\text{\ \ \  (requests)} \\ \text{(LLEs)} \end{matrix}
$ 
\end{tabular}
\label{eq:congestion_problem}
\end{align}
where $\delta, \gamma > 0$ are constants to control the service and admission of requests and LLEs. The following lemma establishes that solving \eqref{eq:congestion_problem} is equivalent to solving \eqref{eq:primal_problem_no_decoherence}. 
\begin{lemma}Suppose $\gamma > \delta > 0$ and that the problem in \eqref{eq:primal_problem_no_decoherence} is feasible. A solution to the problem in \eqref{eq:congestion_problem} is also a solution to the problem in \eqref{eq:primal_problem_no_decoherence}.
\end{lemma}

\begin{IEEEproof}
By construction, the linear objective function enforces to select vectors $z$ and $\hat z$ as large as possible and vectors $x$ and $Ax$ as small as possible. Because $\gamma > \delta$ (admission has preference over service), the optimization will select $z = b \preceq x$ and $\hat z = \hat b = Ax$. That is, the constraints $b \preceq x$ and $Ax \preceq \hat b$ in \eqref{eq:primal_problem_no_decoherence} are satisfied with equality. 
\end{IEEEproof}

\subsection{Algorithm} 
\label{sec:cc_algorithm}
Since \eqref{eq:congestion_problem} is convex, we can use the same approach as in Sec.\ \ref{sec:alg_with_decoherence}: formulate the Lagrange dual problem and solve the static problem with the standard stochastic (sub)gradient ascend method with constant step size (Algorithm \ref{al:cc_spdga}). And as in Sec.\ \ref{sec:alg_with_decoherence}, we can operate the quantum switch if we associate an iteration with a time slot $k=1,2,\dots$ and a state $s_k \in \mathbb S$, where a system state now determines the admission control variables available $B_k := B(s_k)$ and $\hat B_k := \hat B(s_k)$.
The convergence of Algorithm \ref{al:cc_spdga} is standard: it converges to an $\alpha$-neighboring solution (i.e., the optimal average admission/service of requests and LLEs in the switch) where the accuracy can be controlled by making the step size $\alpha$ smaller in Eq.\ \eqref{eq:dual_update_extended} (see, e.g., \cite[Thm.\ 1; claim (i)]{VL19}). However, the convergence of the algorithm is not enough. 
In the next theorem, we show that when $\delta$ and $\gamma$ are set accordingly, Algorithm \ref{al:cc_spdga} has also the desired behavior: it does not attempt to serve requests when there are not enough LLEs available.

\begin{algorithm}[t]
\caption{Stochastic PDGA for the dual of \eqref{eq:congestion_problem}}
\begin{algorithmic}[1]
\State \textbf{Input:} Action set $Y$
\State \textbf{Set:} $\lambda_0 = 0$, $\hat \lambda_0$, $\alpha \in ( 0,1]$
\While{stopping criterion is not met}
\State 1) Observe request and LLE arrivals and construct admission variables sets $B_k$ and $\hat B_k$.
\State 2) Admit requests and LLEs:
\begin{align}
\begin{bmatrix} w_k \\
\hat w_k \end{bmatrix}
\in \underset{v \in B_k, \hat v \in \hat B_k}{\arg \min} \  \left\langle \begin{bmatrix}\lambda_k \\ \hat \lambda_k \end{bmatrix} - \gamma \mathbf 1, \begin{bmatrix} v \\ \hat v \end{bmatrix} \right\rangle \label{eq:admission}
\end{align}
\State 3) Service requests and LLEs: 
\begin{align}
y_k \in \arg \min_{u \in Y} \  \left\langle \delta \mathbf 1- \begin{bmatrix}\lambda_k \\ \hat \lambda_k \end{bmatrix}, \begin{bmatrix} u \\ Au \end{bmatrix} \right\rangle \label{eq:congestion_service}
\end{align}
\State 4) Update the Lagrange multipliers:
\begin{align}
\begin{bmatrix}\lambda_{k+1} \\ \hat \lambda_{k+1} \end{bmatrix}= \left[ \begin{bmatrix}\lambda_k \\ \hat \lambda_k \end{bmatrix} +  \alpha \left(\begin{bmatrix} w_k \\ \hat w_k \end{bmatrix} - \begin{bmatrix} y_k \\ Ay_k \end{bmatrix} \right) \right]^+ \label{eq:dual_update_extended}
\end{align}
\State 5) $k \leftarrow k + 1$
\EndWhile
\end{algorithmic}
\label{al:cc_spdga}
\end{algorithm}

\begin{theorem}
Select $\delta \in [ 2(\gamma + \alpha) + \alpha, 3 (\gamma + \alpha) )$ with $\gamma > 0$ and suppose that the quantum switch can store $ \lceil  \frac{\gamma}{\alpha} + 1 \rceil$ LLEs per client.  Then, Algorithm \ref{al:cc_spdga} will attempt to serve a request from $Q_k$ only if there are enough LLEs in the queues $\hat Q_k$. 
\label{th:main_theorem}
\end{theorem}
\begin{IEEEproof}
To start, note that the Lagrange multipliers in \eqref{eq:dual_update_extended} are equivalent to scaled queues, i.e., $\alpha Q_k = \lambda_k$ and $\alpha \hat Q_k = \hat \lambda_k$. Also, note that $0 \preceq y_k \preceq \mathbf 1$ and $0 \preceq A y_k \preceq \mathbf 1$ (by construction), and that if $y_k \ne 0$, we have that $\langle Ay_k, \mathbf 1\rangle / \langle  y_k, \mathbf 1 \rangle = 2$ since we need two LLEs to serve a request. Eq.\ \eqref{eq:admission} will admit a ``commodity'' if the associated Lagrange multiplier is smaller than or equal to $\gamma$.  Thus, each Lagrange multiplier is upper bounded by $\gamma + \alpha$ for all $k$ and therefore the switch needs to be able to store $\lceil \gamma / \alpha + 1 \rceil$ LLEs per client. 

A sufficient condition for the queues to do not underflow---and so there are always enough LLEs available---is that $\lambda_k - \alpha y_k \succeq 0$ and $\hat \lambda_k - \alpha Ay_k \succeq 0$ for all $k$. We can enforce the latter by selecting $\delta \ge 2(\gamma + \alpha) + \alpha$ in \eqref{eq:congestion_service}, i.e., at least one of the Lagrange multipliers involved in serving a request will have size $\alpha$. Finally, we need that $\delta < 3(\gamma + \alpha)$ since otherwise \eqref{eq:congestion_service} will never attempt to serve a request. 
\end{IEEEproof}

In sum, by selecting $\delta$ and $\gamma$ as indicated in Theorem \ref{th:main_theorem}, Algorithm \ref{al:cc_spdga} has always LLEs available when attempting to serve a request. Finally, note that Algorithm \ref{al:cc_spdga} needs that the quantum switch can store $\lceil \frac{\gamma}{\alpha} + 1 \rceil$ LLEs per client, or conversely, the storage capacity of the switch affects how small we can select $\alpha$ and therefore the algorithm's optimality. Hence, if the switch had infinite memory, we could select $\alpha$ arbitrarily small to maximize the load the switch can support.

\section{Numerical Example}

In this section, we simulate the operation of a quantum switch with $N=6$ clients for the models in Sec.\ \ref{sec:maxdecoherence} and \ref{sec:model_without_decoherence}. 

\subsection{Experiment 1: per time-slot decoherence} \label{simulation:dec}

\begin{figure}
\centering
{\resizebox{0.56\columnwidth}{!}{\input{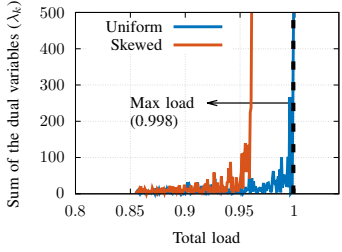}}} 
\caption{Showing the sum of the dual variables at time slot $k = 50000$ for the experiment in Sec.\ \ref{simulation:nodec}.}
\label{fig:1}
\end{figure}

In this experiment, we consider that LLEs are generated successfully in each time slot with probability $0.8$. Also, we consider two types of arrival processes. One where the load is uniformly distributed across requests (Uniform), and another arrival process where the load of three types of requests is 16 times larger than the rest (Skewed).
We run Algorithm \ref{al:spdga} and show in Figure \ref{fig:1} the sum of the dual variables at $k = 50000$ time slots as a function of the total load. Observe from the figure that the sum of the dual variables  is bounded (which is analogous to queue occupancies) and that they grow rapidly as the total arrivals approach the maximum load supported by the switch. This is the expected behavior obtained in queueing systems (e.g., \cite[Sec.\ 4.9]{GNT06}).

\subsection{Experiment 2: no decoherence}
\label{simulation:nodec}

\begin{figure}
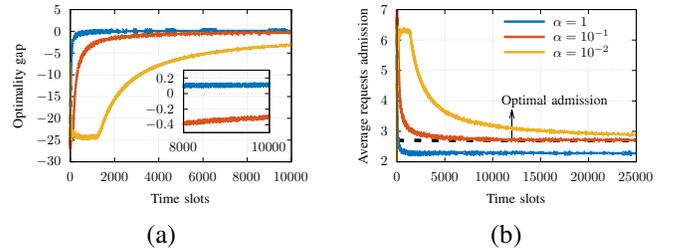

\centering
\begin{tabular}{cc}
{\resizebox{0.47\linewidth}{!}{\input{Figures/exp2.tex}}} & 
{\resizebox{0.47\columnwidth}{!}{\input{Figures/exp3.tex}}} \\
(a) & (b)
\end{tabular}
\caption{Showing the optimality gap and the average admission for the experiment in Sec.\ \ref{simulation:nodec}. The optimal solution corresponds to the solution to the static problem \eqref{eq:congestion_problem}.}
\label{fig:2}
\end{figure}

In this experiment, LLEs and requests arrive with probability $0.3$ in each time slot. We run Algorithm \ref{al:cc_spdga} and show in Figure \ref{fig:2}a the optimality gap as a function of the time slots for different step sizes $\alpha \in \{ 1, 10^{-1}, 10^{-2} \}$. 
Observe from the figure that with $\alpha = 1$,  Algorithm \ref{al:cc_spdga} converges to a nearby optimal point with fewer iterations than with the other step sizes. However, with step sizes $\alpha \in \{ 10^{-1}, 10^{-2} \}$, the optimality gap becomes smaller as the number of time slots increases. The latter can also be observed in Figure \ref{fig:2}b, where we show the average admission of requests in the switch. We can obtain better performance by using smaller step sizes. However, recall that the choice of step size depends on the quantum memory available (Theorem \ref{th:main_theorem} and Sec.\ \ref{sec:cc_algorithm}).


\section{Conclusions}
We have characterized the capacity region of the quantum switch and proposed throughput-optimal policies when LLEs  last (i) for one time slot, or (ii) until they are used to serve a request.  The two models are important as they lower and upper bound the capacity region for any decoherence model. The techniques employed also suggest that a gradient descent based approach might work for decoherence models where LLEs last for an arbitrary time, which is an interesting future research direction.

\section{Acknowledgements}
This work was supported by the European Union’s Horizon 2020 Research and Innovation Program under the Marie Skłodowska-Curie under Agreement 795244. 
The research work was supported by the Army Research Office MURI under the project number W911NF2110325 and by the National Science Foundation under project number EEC-1941583 CQN ERC.

\bibliographystyle{IEEEtran}
\bibliography{references3}

\end{document}